\documentclass[10.5pt,twoside]{article}
\usepackage{amsfonts, amsthm, amssymb,amsmath, txfonts}
\usepackage{latexsym}
\newtheoremstyle{dotless}{}{}{\itshape}{}{\bfseries}{}{ }{}
\theoremstyle{dotless}
\numberwithin{equation}{section}
\newtheorem{thm}{Theorem}[section]

\newtheorem{lemma}{Lemma}[section]

\newtheorem{cor}{Corollary}[section]
\newtheorem{notation}{Notation}


\begin{document}
\title{
\begin{flushleft}
{\small{\it Advanced Nonlinear Studies} {\bf 12} (2012),  ....}
\end{flushleft}
\vspace{0.5in} {\bf\Large  Liouville  Type Theorems for Stable Solutions of\\ Certain Elliptic  Systems}}
\author{{\bf\large Mostafa Fazly\vspace{1mm}} \\
{\it\small Department of Mathematics}\\ {\it\small
University of British Columbia}, {\it\small Vancouver, B.C. Canada V6T 1Z2}\\
{\it\small e-mail: fazly@math.ubc.ca}\vspace{1mm}}
\date{\small{Received 10 December 2010}\\
{\it\small Communicated by  Nassif Ghoussoub}}

\maketitle
\begin{center}
{\bf\small Abstract}

\vspace{3mm} \hspace{.05in}\parbox{4.5in} {{\small We establish Liouville type theorems for elliptic systems with various classes of non-linearities on $\mathbb{R}^N$. We show, among other things, that a system has no semi-stable solution in any dimension,  whenever the infimum of the derivatives of the corresponding non-linearities is positive.
 % in the range of solutions
% then there is no semi-stable solution in any dimension.  After that, we see
We give some immediate applications to various standard systems, such as the Gelfand, and certain Hamiltonian systems. The case where the infimum is zero is more interesting and quite challenging. We show that any  $C^2(\mathbb{R}^N)$ positive entire semi-stable solution of the following Lane-Emden system,
 \begin{eqnarray*}
\hbox{$(N_{\lambda,\gamma})$}\hskip 50pt \left\{ \begin{array}{lcl}
\hfill -\Delta u&=&\lambda  f(x) \ v^p,\\
\hfill -\Delta v&=&\gamma   f(x)  \ u^q,
\end{array}\right.
  \end{eqnarray*}
  is necessarily constant, whenever the  dimension
$N< 8+3\alpha+\frac{8+4\alpha}{q-1}$, provided $p=1$, $q\ge2$ and $f(x)=  (1+|x|^2)^{\frac{\alpha}{2}} $.  The same also holds for $p=q\ge2$ provided $$N < 2+ \frac{2(2+\alpha)}{p-1}  (p+\sqrt{p(p-1)}).$$

 We also consider the case of bounded domains $\Omega\subset\mathbb{R}^N$, where we extend results of  Brown et al. \cite{bs} and Tertikas \cite{te} about stable solutions of equations to systems. At the end, we prove a Pohozaev type theorem for certain weighted elliptic systems.}}
\end{center}


\section{Introduction}

In this paper, we shall examine semi-stable solutions of semi-linear elliptic systems on bounded and unbounded domains on $\mathbb{R}^N$.
%the whole domain and also on the bounded domains.   In fact,
We are particularly interested in Liouville type theorems for stable solutions of such systems. Note that this subject is by now well developed for equations such as those involving Gelfand, Lane-Emden, MEMS and even more general $C^1$ convex non-linearities.
%there are several Liouville type results for semi-linear equations like, Gelfand equation, Lane-Emden equation, MEMS and even with general $C^1$ nonlinearities.
To our knowledge, no such (Liouville type) theorem exists for stable solutions of systems corresponding to these non-linearities. Let us first mention the definition of stability.
% for the most of well-known systems. However, there are several Liouville type results for semilinear equations like, Gelfand equation, Lane-Emden equation, MEMS and even with general $C^1$ nonlinearities.

Let $\Omega$ be a subset of $\mathbb{R}^N$ and $f,g\in C^1(\mathbb{R}^2,\Omega)$. Consider the following general system on $\Omega$,
\begin{eqnarray*}
\label{?}
\hbox{$(N_{f,g})$}\hskip 50pt \left\{ \begin{array}{lcl}
\hfill -\Delta u&=& \ f(u,v,x),   \\
\hfill -\Delta v&=& \ g(u,v,x).
\end{array}\right.
  \end{eqnarray*}
  Following Montenegro \cite{mo}  --whose work was restricted to bounded domains-- we say that a  solution $(u,v)\in C^2(\Omega)$  of  $(N_{f,g})$ is {\em stable} if the following eigenvalue problem, $(S_{f,g})$, has a first positive eigenvalue, $ \eta>0$, with corresponding positive smooth eigenfunction pair $(\phi,\psi)$. We say it is {\em semi-stable} if the first eigenvalue $\eta$ is non-negative, i.e., $\eta\ge0$.
\begin{eqnarray*}
\label{?}
\hbox{$(S_{f,g})$}\hskip 50pt \left\{ \begin{array}{lcl}
\hfill -\Delta \phi &=& \ f_u(u,v,x)\phi  +f_v(u,v,x)\psi +  \eta\phi,   \\
\hfill -\Delta v&=& \ g_u(u,v,x)\phi+g_v(u,v,x)\psi+\eta\psi ,
\end{array}\right.
  \end{eqnarray*}

In section 2, we shall first consider systems on unbounded domains, and establish Liouville type results  for the easy case where the infimum of the derivatives of the non-linearities are greater than zero. We then tackle the more challenging case when the infimum is zero so as to cover Lane-Emden systems.
% we easily can prove Liouville type theorems.
%Then, we see,  by a special case of Lane-Emden system,  if the infimum is zero  proving such theorems are very challenging.

In section 3, we give a general class of non-linearities for which Liouville type theorems can be established.

\section{Systems on unbounded domains}

%In this part,  we give some Liouville type theorems for nonexistence of entire semi-stable solutions in $\mathbb{R}^N$. Let's start with

Assume $\Omega$ is the whole space $\mathbb{R}^N$.  We first prove the following result. %Here is our result for the general system $(N_{f,g})$.
%%%%%%%%%%%%THEOREM 1%%%%5
\begin{thm} \label{unbo}
Assume $f$ and $ g$ are two $C^1$-functions on $\mathbb{R}^{N+2}$.
\begin{enumerate}
\item  If all partial derivatives of $f$ and $ g$ with respect to $u$ and $v$ are non-negative, then there is no $C^2(\mathbb{R}^N)$ entire stable solution for $(N_{f,g})$, i.e., $\eta$ must be zero.
\item If $\min\{\inf f_v,\inf g_u\}\ge C>0$, where the infimum  is on the range of $u$ and $v$, and $f_u,g_v\ge 0$, then there is no $C^2(\mathbb{R}^N)$ entire semi-stable solution for $(N_{f,g})$.
\end{enumerate}
\end{thm}
\noindent{\it Proof.}\, Since proofs of 1 and 2 are similar, we just prove 2. Suppose $(u,v)$ is an entire semi-stable solution of $(N_{f,g})$, so that there exists a pair of positive functions $(\phi,\psi)$ satisfying $(S_{f,g})$ with $\eta\ge 0$.  By adding two equations in $(S_{f,g})$ we get,
%\begin{eqnarray*}
%-\Delta (\phi+\psi)= (f_u+g_u)\phi+(f_v+g_v)\psi+\eta(\phi+\psi).
%\end{eqnarray*}
%If all derivatives of $f$ and $g$ are nonnegative, then
%\begin{eqnarray*}
  % \frac{-\Delta (\phi+\psi)}{\phi+\psi}\ge \eta.
%\end{eqnarray*}
%On the other hand, from the semi-stability assumption $\eta\ge0$, we have
\begin{eqnarray*}
-\Delta (\phi+\psi)= (f_u+g_u)\phi+(f_v+g_v)\psi+\eta(\phi+\psi) \ge  \inf g_u \  \phi+ \inf f_v\ \psi,
\end{eqnarray*}
therefore
\begin{eqnarray*}
\frac{-\Delta (\phi+\psi)}{\phi+\psi}&\ge& (\inf g_u) \frac {\phi} {\phi+\psi }+(\inf f_v)  \frac{\psi} {\phi+\psi }\\
&\ge& \min\{\inf f_v,\inf g_u\}\ge C>0.
\end{eqnarray*}
Now, multiply the above equation by $\zeta^2$ for $0\le\zeta(x)\in C_c^1(\mathbb{R}^N)$, and do integration by parts to get
\begin{eqnarray*}
C \int_{\mathbb{R}^N} \zeta^2  &\le& \int_{\mathbb{R}^N}  \frac{-\Delta (\phi+\psi)}{\phi+\psi}  \zeta^2 \\
&=& \int_{\mathbb{R}^N}   \frac{ 2\zeta}{\phi+\psi}  \nabla(\phi+\psi) \cdot\nabla\zeta -  \frac{ |\nabla(\phi+\psi)|^2}{|\phi+\psi|^2} \zeta^2.
\end{eqnarray*}
Since the function $f(x):=2a x-x^2$ takes it maximum at $x=a$, we conclude
\begin{eqnarray*}
C \int_{\mathbb{R}^N} \zeta^2  \le \int_{\mathbb{R}^N} |\nabla\zeta|^2.
\end{eqnarray*}
This is obviously a contradiction, since there is no Poincar\'e inequality on $\mathbb{R}^N$.
%This can also be seen by considering the following test function,     for $R>0$
%$$\zeta_R(x)=\left\{
 %                     \begin{array}{ll}
    %                    1, & \hbox{if $|x|<R$;} \\
      %                  0, & \hbox{if $|x|>2R$;} \\
 %                          \end{array}
  %                  \right.$$
%in such a way that $||\nabla\zeta_R||_{\infty}\le\frac{c}{R}$.  To get the contradiction, it suffices to take $R>|B_1(0)|(2^N-1)$.

$\hfill \blacksquare$

\vspace{1mm}

\textbf{Examples:} Let $\lambda,\gamma>0$. Consider the following systems on $\mathbb{R}^N$.

\[
\begin{array}{ccc}
\hbox{$(N_{1, \lambda,\gamma})$ }\hskip 20pt \left\{ \begin{array}{lcl}
\hfill -\Delta u &=& \lambda e^v,    \\
\hfill-\Delta v&=& \gamma e^u,
\end{array}\right.
\text{Gelfand System}
\\
\\
 \hbox{$(N_{2, \lambda,\gamma})$ }\hskip 20pt \left\{ \begin{array}{lcl}
-\Delta u&=&\lambda(1+a(x)u+b(x)v)^p, \\
-\Delta v&=&\gamma (1+c(x)u+d(x)v)^q,
\end{array}\right.
 \text{General Lane-Emden System}
\\
\\
\hbox{$(N_{3, \lambda,\gamma})$ }\hskip 20pt  \left\{ \begin{array}{lcl}
\Delta u&=&\lambda v^{-p}, \\
\Delta v&=&\gamma u^{-q},
\end{array}\right.
 \text{General MEMS System}
\\
\\
\hbox{$ (N_{4,\lambda})$ }\hskip 20pt  \left\{ \begin{array}{lcl}
-\Delta u&=& v, \\
-\Delta v&=& \lambda f(u),
\end{array}\right.
\text{Bi-harmonic equation}
\\
\\
 \hbox{$ (N_{5,H}) $ }\hskip 20pt   \left\{ \begin{array}{lcl}
-\Delta u &=& H_v (u,v), \\
-\Delta v &= & H_u (u,v),
\end{array}\right.
 \text{Hamiltonian System}
\end{array}
 \]
where all $C^1(\mathbb{R}^N)$ functions $a,b,c,d,p,q,f, H, H_u$ and $H_v$ are positive.

%%%%%%%%%% need to be finished
\vspace{1mm}

\begin{cor} There is no pair $(u,v)$ of $C^2(\mathbb{R}^N)$ functions satisfying one of the following 4 conditions:
\begin{enumerate}
\item $(u,v)$ is a positive entire semi-stable solution of either $(N_{1,\lambda,\gamma})$ or $(N_{2,\lambda,\gamma})$.
\item $(u,v)$ is a positive bounded entire semi-stable solution of  $(N_{3,\lambda,\gamma})$.
\item $(u,v)$ is an entire semi-stable solution of the bi-harmonic equation  $(N_{4,\lambda}) $ with $f'\ge C>0$ in the range of $u$,
\item $(u,v)$ is an entire semi-stable solution of the Hamiltonian system  $(N_{5,H}) $, where $H_{uv}\ge 0$ and $H_{uu}, H_{vv}\ge C>0$ in the range of $u$ and $v$.
\end{enumerate}
\end{cor}
\newpage
%The conclusion is, if the range of derivatives of nonlinearities $f$ and $g$ in $(N_{f,g})$ are slightly greater than zero, it's not hard to get Liouville type theorems.
 In the next part, we see how much the condition $\min\{\inf f_v,\inf g_u\}\ge C>0$ in Theorem \ref{unbo} on  non-linearities is crucial. Without it, such as the case of Lane-Emden systems, the proof becomes much more delicate and will require restrictions on the dimension $N$.

% We see for a specific example, special case of Lane-Emden system, we need a variety of lemmas to get an appropriate dimension to prove Liouville type theorems.
 %\\Note that there are many papers with nice techniques on stable and unstable solutions of single equations with different nonlinearities which interested readers can find in literature.

\subsection{The Lane-Emden system}

 Existence and non-existence of positive solutions for both Lane-Emden equations and systems have been interesting and challenging questions for decades.

 Recently, Souplet \cite{so} established the Lane-Emden conjecture
in $N=4$ and obtained a new region of non-existence for $N\ge 5$. The so-called Lane-Emden conjecture which has been open for at least fifteen years asserts that the elliptic system
\begin{eqnarray*}
\label{?}
 \left\{ \begin{array}{lcl}
\hfill -\Delta u&=& v^p   \ \ \text{in}\ \ \mathbb{R}^N,\\
\hfill -\Delta v&=&u^q   \ \ \text{in}\ \ \mathbb{R}^N,
\end{array}\right.
  \end{eqnarray*}
for $p,q > 0$ has no positive classical solution if and only if the pair $(p, q)$ lies below the
{\it Sobolev critical hyperbola}, i.e.
$$ \frac{1}{p+1}+\frac{1}{q+1}>1-\frac{2}{N}.$$
This statement is the analogue of the celebrated Gidas-Spruck \cite{gs} Liouville type theorem
for the scalar case. Among other partial results, the conjecture had been proved for $N \le 3$ in \cite{sz, pqs} and for only radial solutions by Mitidieri in \cite{m}.
%Up to now, the conjecture had been proved for radial solutions [5, 8],
%in n  3 space dimensions [9, 6], and in certain subregions below the critical hyperbola
%for n  4 [2, 11, 5, 9, 4, 7, 1].
%I will in particular report on the recent work [10], where I establish the conjecture
%in 4 space dimensions and obtain a new region of nonexistence for n  5.

  Also, the question of the regularity of semi-stable solutions of the system $-\Delta u = \lambda e^v$, $ -\Delta v = \gamma e^u$  in $ \Omega$, and those of the fourth order equation $\Delta^2 u=f(u)$ with zero Dirichlet boundary conditions were examined in \cite{cw} and \cite{cen}, respectively. Roughly speaking, there is a correspondence between the regularity of semi stable solutions on bounded domains and the non-existence of semi stable solutions on $ \mathbb{R}^N$, via rescaling and a blow up procedure.  We combine the techniques of the above mentioned papers to find Liouville type theorems in the notion of stability for a special case of Lane-Emden system in higher dimensions. During preparation of this work, we noticed that Wei-Ye \cite{wy} also used these techniques  to establish independently Liouville type theorems for fourth order equations. To our knowledge, no Liouville type theorem have been established for entire semi-stable solutions of systems of the form
   \begin{eqnarray*}
\label{?}
\hbox{$(N_{\lambda,\gamma})$}\hskip 50pt \left\{ \begin{array}{lcl}
\hfill -\Delta u&=&\lambda  f(x) \ v^p,\\
\hfill -\Delta v&=&\gamma   f(x)  \ u^q,
\end{array}\right.
  \end{eqnarray*}
 where $p\ge1$,  $q>1$ and $\lambda,\gamma\in\mathbb{R^{+}}$.

\begin{notation} {\rm Throughout this subsection, for the sake of simplicity, we say {\it $f \preceq g $} if there exists a positive constant $C$ such that $f\le C g$ holds.}
\end{notation}

We shall first need the following $L^1$-estimates, which were inspired by the work of
% on nonlinearities which are essential and quiet helpful to study the existence and nonexistence of solutions of Lane-Emden system and it's counterparts.
% Interested readers can find more details  in
Serrin and Zou \cite{sz} (whose proof is based on ODE techniques), and Mitidieri and Pohozaev \cite{mp}, who prove similar results by a simpler PDE approach. See also \cite{qs}.
%\\ Here for any  $\lambda,\gamma>0$ and $\alpha\in \mathbb{R}$,  we have an estimate for solutions of  $ (N_{\lambda,\gamma})$.

%%%%%%%%%%%%%%%%%LEMMA 1%%%%%%%%%
\begin{lemma} \label{L1est}
Let $p\ge1, q>1$ and $  \lambda,\gamma>0$. For any  $C^2(\mathbb{R}^N)$ positive entire solution $(u,v)$ of $(N_{\lambda,\gamma})$ and $R>1$, there holds
\begin{eqnarray*}
\int_{B_R}{       f(x)    v^p} &\le&  C_{\lambda,\gamma}  R^{N-\frac{2(q+1)p}{pq-1} - \frac{p+1}{pq-1}\alpha },
\\
\int_{B_R}{    f(x)   u^q}  &\le& C_{\lambda,\gamma}  R^{N-\frac{2(p+1)q}{pq-1} - \frac{q+1}{pq-1}\alpha },
 \end{eqnarray*}
where  $f(x)=  (1+|x|^2)^{\frac{\alpha}{2}}   $ for any $\alpha\in \mathbb{R}$, and $C_{\lambda,\gamma} $ does not depend on $R$.
\end{lemma}

\noindent{\it Proof.}\, Fix the following function $\zeta_R\in C^2_c(\mathbb{R}^N)$ with $0\le\zeta_R\le1$;
 $$\zeta_R(x)=\left\{
                      \begin{array}{ll}
                        1, & \hbox{if $|x|<R$;} \\
                        0, & \hbox{if $|x|>2R$;}
                                                                       \end{array}
                    \right.$$
where $||\nabla\zeta_R||_{\infty}\preceq \frac{1}{R}$ and $||\Delta\zeta_R||_\infty\preceq \frac{1}{R^2}$. For fixed $m\ge 2$, we have

$$|\Delta \zeta^m_R(x)|\preceq \left\{
                      \begin{array}{ll}
                        0, & \hbox{if $|x|<R$ or $|x|>2R$;} \\
                         R^{-2} \zeta^{m-2}_R , & \hbox{if $R<|x|<2R$.}
                                                                       \end{array}
                    \right.$$
For $m\ge 2$, test the first equation of  $(N_{\lambda,\gamma})$ by $\zeta^m_R$ and integrate to get
\begin{eqnarray*}
\lambda \int_{\mathbb{R}^N}     (1+|x|^2)^{\frac{\alpha}{2}}       v^p  \zeta^m_R &= &- \int_{\mathbb{R}^N}  \Delta u   \zeta^m_R\\
&=& - \int_{\mathbb{R}^N} u\Delta \zeta^m_R   \preceq  R^{-2}\int_{B_{2R}\setminus B_{R}} u\zeta^{m-2}_R.
\end{eqnarray*}
Applying H\"{o}lder's inequality we get
\begin{eqnarray*}
\lambda \int_{\mathbb{R}^N}     (1+|x|^2)^{\frac{\alpha}{2}}       v^p  \zeta^m_R &\preceq &  R^{-2}  \left(         \int_{B_{2R}\setminus B_{R}}  (1+|x|^2)^{\frac{-\alpha}{2q }q'    }          \right)^{\frac{1}{q'}}    \\
&&\left(     \int_{B_{2R}\setminus B_{R}}      (1+|x|^2)^{\frac{\alpha}{2}}      u^q \zeta^{(m-2)q}_R   \right)^{1/q}\\
&\preceq &  R^{  (N-\frac{\alpha}{q}q')\frac{1}{q'}  -2 }  \\
&& \left(     \int_{B_{2R}\setminus B_{R}}      (1+|x|^2)^{\frac{\alpha}{2}}      u^q \zeta^{(m-2)q}_R   \right)^{1/q}.
\end{eqnarray*}
By a similar calculation for $k\ge 2$, we obtain
\begin{eqnarray*}
\gamma \int_{\mathbb{R}^N}     (1+|x|^2)^{\frac{\alpha}{2}}       u^q  \zeta^k_R   &\preceq &  R^{-2}\int_{B_{2R}\setminus B_{R}} v\zeta^{k-2}_R     (1+|x|^2)^{\frac{\alpha}{2p}}    (1+|x|^2)^{\frac{-\alpha}{2p}}.
\end{eqnarray*}
 %If $p=1$, we have
%\begin{eqnarray*}
%\gamma \int_{\mathbb{R}^N}     (1+|x|^2)^{\frac{\alpha}{2}}       u^q  \zeta^k_R &\le& c %R^{-2-\alpha} \int_{B_{2R}\setminus B_{R}} v\zeta^{k-2}_R     (1+|x|^2)^{\frac{\alpha}{2}}  ,
% \end{eqnarray*}
%and for $p>1$
By H\"{o}lder's inequality we get
\begin{eqnarray*}
\gamma \int_{\mathbb{R}^N}     (1+|x|^2)^{\frac{\alpha}{2}}       u^q  \zeta^k_R &\preceq &R^{  (N-\frac{\alpha}{p}p')\frac{1}{p'}  -2 }     \left(\int_{B_{2R}\setminus B_{R}}    (1+|x|^2)^{\frac{\alpha}{2}} v^p \zeta^{(k-2)p}_R   \right)^\frac{1}{p} ,
 \end{eqnarray*}
where $\frac{1}{p}+\frac{1}{p'}=1$. For $p=1$ take $p'=\infty$.
 For large enough $k$ we have $2+\frac{k}{q}<(k-2)p$, so we can choose $m$ such that $2+\frac{k}{q}\le   m  \le   (k-2)p$ which says $m\le (k-2)p$ and $k\le (m-2)q$. Therefore, by collecting  above inequalities we get for $p\ge 1$,
  %\begin{eqnarray*}
%\gamma \int_{\mathbb{R}^N}     (1+|x|^2)^{\frac{\alpha}{2}}       u^q  \zeta^k_R &\le& c %R^{-2-\alpha}   R^{  (N-\frac{\alpha}{q}q')\frac{1}{q'}  -2 } \left(     \int_{B_{2R}\setminus B_{R}}      %(1+|x|^2)^{\frac{\alpha}{2}}      u^q \zeta^{k}_R   \right)^{1/q},
%\\
%\lambda \int_{\mathbb{R}^N}     (1+|x|^2)^{\frac{\alpha}{2}}       v  \zeta^m_R &\le& c %R^\frac{-2-\alpha}{q}   R^{  (N-\frac{\alpha}{q}q')\frac{1}{q'}  -2 } \left(     \int_{B_{2R}\setminus %B_{R}}      (1+|x|^2)^{\frac{\alpha}{2}}      v \zeta^{m}_R   \right)^{1/q},
%\end{eqnarray*}
%and, for $p>1$
  \begin{eqnarray*}
\gamma \int_{\mathbb{R}^N}     (1+|x|^2)^{\frac{\alpha}{2}}       u^q  \zeta^k_R &\preceq & R^{  (N-\frac{\alpha}{p}p')\frac{1}{p'}  -2 } R^{  \left(  (N-\frac{\alpha}{q}q')\frac{1}{q'}  -2 \right)\frac{1}{p}}  \\
&& \left(     \int_{B_{2R}\setminus B_{R}}      (1+|x|^2)^{\frac{\alpha}{2}}      u^q \zeta^{k}_R   \right)^{1/pq},
\\
\lambda \int_{\mathbb{R}^N}     (1+|x|^2)^{\frac{\alpha}{2}}       v^p  \zeta^m_R &\preceq &  R^{ \left( (N-\frac{\alpha}{p}p')\frac{1}{p'}  -2  \right)\frac{1}{q} }   R^{  (N-\frac{\alpha}{q}q')\frac{1}{q'}  -2 } \\
&&\left(     \int_{B_{2R}\setminus B_{R}}      (1+|x|^2)^{\frac{\alpha}{2}}      v^p \zeta^{m}_R   \right)^{1/pq}.
\end{eqnarray*}
Therefore,  for $p\ge1$ we have
 \begin{eqnarray*}
\int_{B_R}{         (1+|x|^2)^{\frac{\alpha}{2}}       v^p} &\preceq &    R^{N-\frac{2(q+1)p}{pq-1} - \frac{p+1}{pq-1}\alpha },
\\
\int_{B_R}{    (1+|x|^2)^{\frac{\alpha}{2}}   u^q}  &\preceq &   R^{N-\frac{2(p+1)q}{pq-1} - \frac{q+1}{pq-1}\alpha }.
 \end{eqnarray*}

$\hfill \blacksquare$ %\hfill $\Box$

%%%%%%%%%%%%%%%Corollary 2%%%%%%%%%
\begin{cor} \label{cor}
With the same assumptions as Lemma \ref{L1est}, we have
 \begin{eqnarray*}
\int_{B_R}{        \frac{u^2}{ (1+|x|^2)^{\frac{\alpha}{2}}  }     } &\le&  C  R^{N-\frac{4(p+1)}{pq-1} - \frac{2(p+1)}{pq-1}\alpha-\alpha },
 \end{eqnarray*}
where the positive constant $C$ does not depend on $R$.
\end{cor}

\noindent{\it Proof.} Apply H\"{o}lder's inequality  to obtain
\begin{eqnarray*}
\int_{B_R}{          \frac{u^2}{ (1+|x|^2)^{\frac{\alpha}{2}}   } } &\le&  \left ( \int_{B_R}    (1+|x|^2)^{\frac{\alpha}{2}}       u^q  \right)^{ \frac{2}{q} }   \left(   \int_{B_R}   (1+|x|^2)^{ (-\frac{\alpha}{2}   - \frac{\alpha}{q}  )\frac{q}{q-2} }  \right)^{1-\frac{2}{q} }\\
& \preceq &  R^{N(1-\frac{2}{q})-(1+\frac{2}{q})\alpha} \left( \int_{B_R} (1+|x|^2)^{\frac{\alpha}{2}} u^q   \right)^{\frac{2}{q}}.
\end{eqnarray*}
Now, use Lemma \ref{L1est} to get the desired inequality.

$\hfill \blacksquare$ %\hfill $\Box$

An immediate application of Lemma \ref{L1est} is the following Liouville type theorem, in the absence of stability, for $(N_{\lambda,\gamma})$ which is related to the Lane-Emden conjecture.  In the case $\alpha=0$, this was obtained by Mitidieri in \cite{m} and also modified by Serrin and Zou in \cite{sz}.

\vspace{1mm}
\begin{thm} \label{louvill}
Let $p\ge1, q>1$ and $  \lambda,\gamma>0$. Assume $(u,v) $ is a non-negative entire classical solution for $(N_{\lambda,\gamma})$, provided
 \begin{eqnarray*}
N-2 \le \max\left\{  \frac{(2+\alpha)(q+1)}{pq-1} , \frac{(2+\alpha)(p+1)}{pq-1}\right\}.
 \end{eqnarray*}
Then, $(u,v)$ must be the trivial solution.
\end{thm}

\noindent{\it Proof.} Proof is a direct consequnce of Lemma \ref{L1est}.

$\hfill \blacksquare$ %\hfill $\Box$

%%%%%%%%%%%%%%%%Lemma 2%%%%%%%%
Now, we prove a point-wise comparison result for solutions of $(N_{\lambda,\gamma})$. The following lemma is an adaptation of a result  of Souplet \cite{so}.
\begin{lemma} \label{souplet}
  \, If $(u,v)$ is a $C^2(\mathbb{R}^N)$  positive entire solution of $(N_{\lambda,\gamma})$, then we have
 \begin{eqnarray*}
u^{q+1}\le \frac{q+1}{p+1}v^{p+1},
\end{eqnarray*}
where $q\ge \max\{p,2\}$, $\lambda\le\gamma$ and $\alpha>-2$.
\end{lemma}
\noindent{\it Proof.}\,  Define $w:=u-\beta v^t$ for $0<t\le1$ and $\beta>0$, so we have
\begin{eqnarray*}
 \Delta w &=&\Delta u - \beta  t \Delta v \ v^{t-1}-\beta  t(t-1)|\nabla v|^2v^{t-2}\\
 &\ge& -\lambda    (1+|x|^2)^{\frac{\alpha}{2}}     v^p+\beta \ t \gamma  (1+|x|^2)^{\frac{\alpha}{2}}  u^qv^{t-1}\\
 & \ge& \gamma  (1+|x|^2)^{\frac{\alpha}{2}} (-v^p+\beta ^{q+1} t\ v^{tq+t-1}), \ \ \ \ \ \ \text{on }  \{w\ge 0 \}.
 \end{eqnarray*}
By taking $t=\frac{p+1}{q+1}$ and $\beta \ge (\frac{q+1}{p+1})^{\frac{1}{q+1}}$, we have $\Delta w\ge 0$ on $\{w\ge 0\}$. Now, by Green's theorem we have
\begin{eqnarray}  \label{w}
\int_{B_R}{|\nabla w_+|^2}&=&-\int_{B_R}{w_+\Delta w} + \int_{|x|=R}{w_+(x)\partial_\nu w(x) \ dS(x)}\\
&\le& R^{N-1}\int_{|z|=1}{w_+(Rz)w_r(Rz) dS(z)}=\frac{1}{2}R^{N-1}g'(R),
\end{eqnarray}
 where
$g(R)=\int_{|z|=1}{w_+^2(Rz) dS(z)}$. Moreover, on the set $\{w\ge0\}$ we have $\beta v^t\le u$ and $u$ and $v$ are positive, so
\begin{eqnarray*}
g(R)&=&\int_{ \{|z|=1\},\{w\ge0\} } {  (u(Rz)-\beta v^t(Rz)    )^2\ dS(z)    }\\
& \le& \int_{  \{|z|=1\}   }{  u^2(Rz)\ dS(z)  }\\
&\preceq & \left(\int_{  \{|z|=1\}   }{  u^q(Rz)\ dS(z)  }\right)^{2/q},\  \ \ \ \ \text{since } \ \ q\ge2.
\end{eqnarray*}
Also, by Lemma \ref{L1est}, we know that
 \begin{eqnarray*}
 \int_{0}^{R}r^{N-1}\int_{  \{|z|=1\}   }{     (1+|rz|^2)^{\frac{\alpha}{2}}     u^q(rz)\ dS(z) dr } &=&     \int_{0}^{R}r^{N-1}   (1+r^2)^{\frac{\alpha}{2}} \\
 &&
 \int_{  \{|z|=1\}   }{u^q(rz)\ dS(z) dr }     \\
 &\le& C_{\lambda,\gamma}  R^{N-\frac{2(p+1)q}{pq-1}- \frac{q+1}{pq-1}\alpha      } .
 \end{eqnarray*}
Therefore, for $\alpha>-2$,  there is a sequence $R_i\to \infty$ such that $$\int_{  \{|z|=1\}   }{  u^q(R_iz)\ dS(z)  }\to 0.$$
This means  that $g(R_i)\to 0$ for $R_i\to \infty$. Since $g$ is a positive function, there is a sequence $\hat{R}_i\to \infty$ such that $g'(\hat{R}_i)\le0$.  Hence, from (\ref{w}) we see that $\int_{B_{\hat{R}_i}}{|\nabla w_+|^2}=0$, i.e. $w_+$ is constant. If $w_+
=C>0$, then by continuity we conclude $w=C$ and by definition of $w$ we see $u\ge C>0$ in $\mathbb{R}^N$ which is in contradiction with Lemma 1. Therefore, $w_+=0$ and $u\le \beta v^t$.

$\hfill \blacksquare$ %\hfill $\Box$

\begin{cor} Let $(u,v)$ be a $C^2(\mathbb{R}^N)$ entire positive solution  of  $(N_{\lambda,\gamma})$ with $p=q\ge 2$, then $u=(\frac{\lambda}{\gamma})^{\frac{1}{p+1}}  v$.
 \end{cor}

It follows that $(N_{\lambda,\gamma})$ with  $p=q\ge2$ reduces to the single equation which has been classified by Farina in \cite{f} for $\alpha=0$. Similar ideas as in Farina's easily yield the following result.

%%%%%%%%%%%%%%THEOREM 2%%%%%%%%
\begin{thm}
If $(u,v)$ is a $C^2(\mathbb{R}^N)$ non-negative entire semi-stable solution of $  (N_{\lambda,\gamma})$ with $p=q\ge2$, and
$$N < 2+ \frac{2(2+\alpha)}{p-1}  (p+\sqrt{p(p-1)}), $$
then, $(u,v)$ is the trivial solution.
\end{thm}

Similar results were also obtained by  Esposito-Ghoussoub-Guo \cite{eggu1, eggu2} %for nonlinearities with negative exponents
and by Esposito \cite{e1}. They can be applied to obtain the following analogue for systems.

%%%%%%%%%%%%%%THEOREM 3%%%%%%%%

%Here is our main result of this section. The idea of proof has been derived from Cowan-Esposito-Ghoussoub \cite{cen} and Wei-Ye \cite{wy}.
\begin{thm} \label{main}
Suppose $(u,v)$ is a $C^2(\mathbb{R}^N)$ non-negative entire semi-stable solution of  $(N_{\lambda,\gamma})$ with $p=1$ in dimensions
\begin {equation} \label {assum}
N< 8+3\alpha+\frac{8+4\alpha}{q-1}.
\end{equation}
 Then, $(u,v)$ is the trivial solution.
\end{thm}

%%%%%%%%%%%%%%%%Lemma 3%%%%%%%%

First, a simple calculation leads us to the following identity.
\begin{lemma} \label{identity}
For any $w\in C^2(\mathbb{R}^N)$ and $\zeta\in C_c^2(\mathbb{R}^N)$, we have
\begin{eqnarray*}
\frac{\Delta w}{   (1+|x|^2)^{\frac{\alpha}{2}}  } \Delta (w\zeta^2)-\frac{ |\Delta (w\zeta)|^2 }{  (1+|x|^2)^{\frac{\alpha}{2}}  }&=&-4 \frac{| \nabla w\cdot\nabla\zeta  |^2}{  (1+|x|^2)^{\frac{\alpha}{2}}   } -\frac{w^2}{  (1+|x|^2)^{\frac{\alpha}{2}}  } |\Delta \zeta|^2\\
&+& 2 \frac{w\Delta w}{ (1+|x|^2)^{\frac{\alpha}{2}} } |\nabla \zeta|^2 -2 \frac{\nabla w^2 \cdot\nabla\zeta   }{ (1+|x|^2)^{\frac{\alpha}{2}} } \Delta \zeta.
\end{eqnarray*}
\end{lemma}

\noindent{\it Proof.}\,
Let $\lambda=\gamma=1$. Assume $(u,v)$ is a semi-stable positive solution. Inspired by the linearised equation
\begin{eqnarray*}
-\Delta \phi&=& (1+|x|^2)^{\frac{\alpha}{2}} \psi,\\
-\Delta \psi&=& q(1+|x|^2)^{\frac{\alpha}{2}} u^{q-1} \phi,
\end{eqnarray*}
 we have the following stability inequality, for all $\zeta\in C_c^2(\mathbb{R}^N)$
\begin{eqnarray} \label{stability}
q\int_{\mathbb{R}^N }    (1+|x|^2)^{\frac{\alpha}{2}}   u^{q-1}\zeta^2 \le \int_{\mathbb{R}^N}{    \frac{|\Delta\zeta|^2}{    (1+|x|^2)^{\frac{\alpha}{2}}     }}.
\end{eqnarray}

\textbf{Step 1. Test the stability inequality on $u$}. Set $\zeta=u\xi$ for $\xi\in  C_c^2(\mathbb{R}^N)$, and test the stability inequality on $\zeta $ to get
\begin{eqnarray*}
(q-1) \int_{\mathbb{R}^N }    (1+|x|^2)^{\frac{\alpha}{2}}   u^{q-1}|u\xi|^2 &\le&  \int_{\mathbb{R}^N}{    \frac{|\Delta(u\xi)|^2}{    (1+|x|^2)^{\frac{\alpha}{2}}     }      -      (1+|x|^2)^{\frac{\alpha}{2}}   u^q  u\xi^2}\\
&=& \int_{\mathbb{R}^N}    \frac{|\Delta(u\xi)|^2}{    (1+|x|^2)^{\frac{\alpha}{2}}     }
+\int_{\mathbb{R}^N}  u\Delta v \xi^2\\
&=& \int_{\mathbb{R}^N}    \frac{|\Delta(u\xi)|^2}{    (1+|x|^2)^{\frac{\alpha}{2}}     }
+\int_{\mathbb{R}^N}  v\Delta  (u\xi^2)\\
&=&\int_{\mathbb{R}^N}    \frac{|\Delta(u\xi)|^2}{    (1+|x|^2)^{\frac{\alpha}{2}}     }
-\int_{\mathbb{R}^N}  \frac{\Delta u}{  (1+|x|^2)^{\frac{\alpha}{2}}    }              \Delta (u \xi^2).
\end{eqnarray*}
Now, using Lemma \ref {identity}, we get
\begin{eqnarray*}
(q-1) \int_{\mathbb{R}^N }    (1+|x|^2)^{\frac{\alpha}{2}}   u^{q-1}|u\xi|^2 &\le&  4 \int_{\mathbb{R}^N} \frac{  |\nabla u|^2 |\nabla \xi|^2}{   (1+|x|^2)^{\frac{\alpha}{2}}  } +  \int_{\mathbb{R}^N } \frac{  u^2}{   (1+|x|^2)^{\frac{\alpha}{2}}  }  |\Delta\xi|^2\\
&+& 2\int_{\mathbb{R}^N } uv |\nabla\xi|^2 +2  \int_{\mathbb{R}^N } u^2 div (\frac{\nabla\xi\Delta \xi}{   (1+|x|^2)^{\frac{\alpha}{2}}   }).
\end{eqnarray*}
By Green's theorem, we can modify the first term in the right hand side as
\begin{eqnarray*}
 \int_{\mathbb{R}^N} \frac{  |\nabla u|^2 |\nabla \xi|^2}{   (1+|x|^2)^{\frac{\alpha}{2}}  }&=&  \int_{\mathbb{R}^N}   \frac{u(-\Delta u)}{     (1+|x|^2)^{\frac{\alpha}{2}}   }   |\nabla\xi|^2         \\
 &+& \frac{1}{2}\int_{\mathbb{R}^N} u^2    div \left(\frac{\nabla|\nabla\xi|^2}{   (1+|x|^2)^{\frac{\alpha}{2}}  }     +\nabla (1+|x|^2)^{\frac{-\alpha}{2}} |\nabla\xi|^2   \right),
 \end{eqnarray*}
combine this equality with the last inequality to get
\begin{eqnarray*}
 \int_{\mathbb{R}^N }    (1+|x|^2)^{\frac{\alpha}{2}}   u^{q-1}|u\xi|^2 &\preceq &  \int_{\mathbb{R}^N } uv|\nabla\xi|^2 \\
 &+& \int_{\mathbb{R}^N } \frac{u^2}{   (1+|x|^2)^{\frac{\alpha}{2}}  } \left( |\Delta\xi|^2+|\nabla\xi\cdot\nabla\Delta\xi|+|\Delta|\nabla\xi|^2|\right)\\
&+& \int_{\mathbb{R}^N }   \frac{u^2}{   (1+|x|^2)^{\frac{\alpha+1}{2}}  }  \left( |\nabla\xi||\Delta\xi|+|\nabla|\nabla\xi|^2|\right)\\
&+& \int_{\mathbb{R}^N } \frac{u^2}{   (1+|x|^2)^{\frac{\alpha+2}{2}}  }   |\nabla\xi|^2\\
&=:&  \int_{\mathbb{R}^N } uv|\nabla\xi|^2 +I(u,\xi).
 \end{eqnarray*}
 From this and the stability inequality, (\ref{stability}), we conclude
\begin{eqnarray}\label{ine3}
 \int_{\mathbb{R}^N }   \frac{|\Delta(u\xi)|^2}{ (1+|x|^2)^{\frac{\alpha}{2}}  }&\preceq &    \int_{\mathbb{R}^N } uv|\nabla\xi|^2 + I(u,\xi).
 \end{eqnarray}
Using the following equality
\begin{eqnarray*}
 \Delta(u\xi)=- (1+|x|^2)^{\frac{\alpha}{2}} v\xi+u\Delta\xi+2\nabla u\cdot\nabla\xi,
 \end{eqnarray*}
 and (\ref{ine3}) we get
\begin{eqnarray*}
     \int_{\mathbb{R}^N }   (1+|x|^2)^{\frac{\alpha}{2}}  v^2 \xi^2 &\preceq &    \int_{\mathbb{R}^N } uv|\nabla\xi|^2 + I(u,\xi).
 \end{eqnarray*}
Now, set $\xi=\zeta^m_R\in C^2_c(\mathbb{R}^N)$ for $0\le\zeta_R\le1$,  $m>2$,  and
 $$\zeta_R(x):=\left\{
                      \begin{array}{ll}
                        1, & \hbox{if $|x|<R$;} \\
                        0, & \hbox{if $|x|>2R$;}
                                                                       \end{array}
                    \right.$$
where $||\nabla\zeta_R||_{\infty}\preceq \frac{1}{R}$ and $||\Delta\zeta_R||_\infty \preceq \frac{1}{R^2}$. So,  by Young's inequality we get
\begin{eqnarray*}
 \int_{\mathbb{R}^N }   (1+|x|^2)^{\frac{\alpha}{2}}  v^2 \zeta_R^{2m} &\preceq &    \int_{B_{2R}\setminus B_{R} } uv|\nabla\zeta_R|^2 \zeta_R^{2(m-1)} + I(u,\zeta^m_R)\\
 &\le& \frac{1}{2C}  \int_{B_{2R}\setminus B_{R} }   (1+|x|^2)^{\frac{\alpha}{2}}  v^2 \zeta_R^{2m}
 \\
 &+&C  \int_{B_{2R}\setminus B_{R} }  \frac{u^2 }{(1+|x|^2)^{\frac{\alpha}{2}} }  |\nabla\zeta_R|^4 \zeta_R^{2(m-2)}+   I(u,\zeta^m_R).
 \end{eqnarray*}
Therefore,
\begin{eqnarray} \label{ine1}
 \int_{\mathbb{R}^N }   (1+|x|^2)^{\frac{\alpha}{2}}  v^2 \zeta_R^{2m} &\preceq &  R^{-4} \int_{B_{2R}\setminus B_{R} } \frac{u^2 }{(1+|x|^2)^{\frac{\alpha}{2}} } \zeta_R^{2(m-2)}+  I(u,\zeta^m_R).
      \end{eqnarray}
On the other hand, by definition of $I(.,.)$ we have
\begin{eqnarray*}
I(u,\zeta^m_R) &\le&     R^{-4}   \int_{B_{2R}\setminus B_{R} } \frac{u^2 }{(1+|x|^2)^{\frac{\alpha}{2}} } \zeta_R^{2(m-2)} \\
&+&R^{-3}  \int_{B_{2R}\setminus B_{R} } \frac{u^2 }{(1+|x|^2)^{\frac{\alpha+1}{2}} } \zeta_R^{2(m-2)}  \\
&+&R^{-2} \int_{B_{2R}\setminus B_{R} } \frac{u^2 }{(1+|x|^2)^{\frac{\alpha+2}{2}} } \zeta_R^{2(m-2)}  \\
&\preceq &  R^{-4}   \int_{B_{2R}\setminus B_{R} } \frac{u^2 }{(1+|x|^2)^{\frac{\alpha}{2}} } \zeta_R^{2(m-2)}.
\end{eqnarray*}
 From this and  (\ref{ine1}), we get
\begin{eqnarray*}
 \int_{\mathbb{R}^N }   (1+|x|^2)^{\frac{\alpha}{2}}  v^2 \zeta_R^{2m} &\preceq &   R^{-4} \int_{B_{2R}\setminus B_{R} } \frac{u^2 }{(1+|x|^2)^{\frac{\alpha}{2}} } \zeta_R^{2(m-2)},
  \end{eqnarray*}
in the light of Corollary \ref{cor}, we see
\begin{eqnarray} \label{L2estimate}
 \int_{\mathbb{R}^N }   (1+|x|^2)^{\frac{\alpha}{2}}  v^2 \zeta_R^{2m} &\preceq &   R^{N-4-\alpha-\frac{8+4\alpha}{q-1}     }.
  \end{eqnarray}

\textbf{Step 2. Test the stability inequality on $v$}.   Set $\zeta=v\xi$ for $\xi\in  C_c^2(\mathbb{R}^N)$ and test the stability inequality,  (\ref{stability}),  on $\zeta $ to get
\begin{eqnarray*}
0&\le&  \int_{\mathbb{R}^N }    \frac{|\Delta (v\xi)|^2}{    (1+|x|^2)^{\frac{\alpha}{2}}  } - q  (1+|x|^2)^{\frac{\alpha}{2}}   u^{q-1} |v\xi|^2\\
 &=&   \int_{\mathbb{R}^N }    \frac{|\Delta (v\xi)|^2}{(1+|x|^2)^{\frac{\alpha}{2}}  } + q    u^{q-1} \Delta u \ v \xi^2\\
&=&  \int_{\mathbb{R}^N }    \frac{|\Delta (v\xi)|^2}{(1+|x|^2)^{\frac{\alpha}{2}}  } +\Delta (u^q)v\xi^2-q(q-1)|\nabla u|^2 u^{q-2} v\xi^2.
\end{eqnarray*}
Therefore,
\begin{eqnarray*}
q(q-1) \int_{\mathbb{R}^N }   |\nabla u|^2 u^{q-2} v\xi^2 &\le&   \int_{\mathbb{R}^N }    \frac{|\Delta (v\xi)|^2}{(1+|x|^2)^{\frac{\alpha}{2}}  } +   \int_{\mathbb{R}^N }  \Delta (u^q)v\xi^2\\
& =&\int_{\mathbb{R}^N }    \frac{|\Delta (v\xi)|^2}{(1+|x|^2)^{\frac{\alpha}{2}}  } +   \int_{\mathbb{R}^N }  u^q\Delta (v\xi^2)\\
&= &\int_{\mathbb{R}^N }    \frac{|\Delta (v\xi)|^2}{(1+|x|^2)^{\frac{\alpha}{2}}  } -  \int_{\mathbb{R}^N }  \frac{\Delta v}{   (1+|x|^2)^{\frac{\alpha}{2}}     }   \Delta(v\xi^2).
 \end{eqnarray*}
By the same idea as Step 1 and using Lemma \ref{identity} we get
\begin{eqnarray}\label{ine2}
 \int_{\mathbb{R}^N }   |\nabla u|^2 u^{q-2} v\xi^2 &\preceq &   \int_{\mathbb{R}^N } u^q v |\nabla\xi|^2 +  I(v,\xi).
 \end{eqnarray}
On the other hand, by following the ideas of Cowan-Esposito-Ghoussoub \cite{cen}, we have
\begin{eqnarray*}
\int_{\mathbb{R}^N }  h(u)(-\Delta u) \xi^2= \int_{\mathbb{R}^N } h'(u) |\nabla u|^2 \xi^2 -\int_{\mathbb{R}^N } H(u) \Delta \xi^2,
 \end{eqnarray*}
where $H(u)=\int_{0}^{u} h(t) dt$. Let $h(u):=u^{\frac{3q-1}{2}}$,  use Lemma \ref{souplet} to get
\begin{eqnarray*}
   \int_{\mathbb{R}^N }      (1+|x|^2)^{\frac{\alpha}{2}}    u^{2q} \xi^2 &\preceq &   \int_{\mathbb{R}^N }    (1+|x|^2)^{\frac{\alpha}{2}}  u^{\frac{3q-1}{2}} v \xi^2   \\
   &\preceq &  \int_{\mathbb{R}^N }   u^{\frac{3q-3}{2}} |\nabla u|^2 \xi^2+ \int_{\mathbb{R}^N }  u^{\frac{3q+1}{2}} |\Delta \xi^2| \\
   &\preceq &  \int_{\mathbb{R}^N }   u^{q-2} v|\nabla u|^2 \xi^2+  \int_{\mathbb{R}^N }  u^{q} v |\Delta \xi^2|.
 \end{eqnarray*}
From this and (\ref{ine2}) we see that
\begin{eqnarray*}
   \int_{\mathbb{R}^N }      (1+|x|^2)^{\frac{\alpha}{2}}    u^{2q} \xi^2 & \preceq &   \int_{\mathbb{R}^N }  u^{q} v(|\nabla\xi|^2+ |\Delta \xi^2|) +  I(v,\xi).
 \end{eqnarray*}
Using the same test function as Step 1 and doing similar calculation we get
\begin{eqnarray*}
   \int_{\mathbb{R}^N }      (1+|x|^2)^{\frac{\alpha}{2}}    u^{2q} \zeta_R^{2m} & \preceq & R^{-4}  \int_{B_{2R}\setminus B_{R} }  \frac{v^2}{    (1+|x|^2)^{\frac{\alpha}{2}}      } \zeta_R^{2(m-2)}\\
   &\preceq &  R^{-4-2\alpha}  \int_{B_{2R}}    (1+|x|^2)^{\frac{\alpha}{2}} v^2  .
    \end{eqnarray*}
By the result of Step 1, (\ref{L2estimate}), we get
\begin{eqnarray*}
   \int_{B_{2R} }      (1+|x|^2)^{\frac{\alpha}{2}}    u^{2q} & \preceq &  R^{N-8-3\alpha-\frac{8+4\alpha}{q-1}       }  .
     \end{eqnarray*}

$\hfill \blacksquare$ %\hfill $\Box$

\noindent{\bf Remark}  One can see that Lemma \ref{L1est} and \ref{souplet} can be adopted for the following system.  So, a counterpart of Theorem \ref{main} can be proved.
 \begin{eqnarray*}
\label{?}
\hbox{$(N'_{\lambda,\gamma})$}\hskip 50pt \left\{ \begin{array}{lcl}
\hfill -\Delta u&=&\lambda  f(x) \ v,\\
\hfill -\Delta v&=&\gamma   f(x)  \ u^q + \gamma f(x) \ u^r.
\end{array}\right.
  \end{eqnarray*}

\noindent{\bf Open problem}:  A natural question is how we can establish a Liouville type theorem for $(N_{\lambda,\gamma})$ with other values of $p$ and $q$.

\section{Systems on  bounded domains}
In this part, we generalize results of Brown et al. \cite{bs} and Tertikas \cite{te}  about  the stable solutions of equations  to systems.   Consider the following system:

\begin{eqnarray*}
\label{?}
\hbox{$(N_{f,g})$}\hskip 50pt \left\{ \begin{array}{lcl}
\hfill -div(e^\theta \nabla u) &=&\lambda e^\theta f(v)  \ \ \ \text{in}  \ \Omega,  \\
\hfill -div(e^\theta \nabla v) &=&\gamma e^\theta g(u)   \  \ \ \text{in}  \ \Omega, \\
\end{array}\right.
    \end{eqnarray*}
with the Robin  boundary conditions:
  \begin{eqnarray*}
 \left\{ \begin{array}{lcl}
\hfill au+b\partial_\nu u &=&0  \ \  \text{on} \  \partial\Omega, \\
\hfill  av+b\partial_\nu v &=&0 \ \ \text{on} \  \partial\Omega, \\
\end{array}\right.
  \end{eqnarray*}
%\[
%\begin{array}{ccc}
%-div(e^\theta \nabla u) =\lambda e^\theta f(v)   &  \text{in}  \ \Omega &  \ \ \ \  (N_{\lambda,\gamma})\\
%-div(e^\theta \nabla v) =\gamma e^\theta g(u)   &  \text{in} \ \Omega &   \\
% &\\
 %au+b\partial_\nu u=0 &  \text{on} \  \partial\Omega  &   \\
  %av+b\partial_\nu v=0 &  \text{on} \  \partial\Omega  &
%\end{array}
% \]
where $\Omega$ is an open bounded subset of $\mathbb{R}^N$,  $f,g\in C^2(\mathbb{R})$,  $\theta\in C^1(\Omega)$,    $a,b\in\mathbb{R}$ and $\lambda,\gamma\in\mathbb{R}^+$.   In the special case $\theta=0$, $\lambda=\gamma$ and $f=g$ we have the single equation.

We define $\mathbb{X}$ to be the following  subset of convex functions in $C^2(\mathbb{R})$: $h\in\mathbb{X}$ if we have
\[
\begin{array}{ccc}
\text{either} &\ h(0)<0  \ \text{and }\  h''(u)\ge 0 \ \ \text{for} \  u>0,&   \\
\text{or} &\ h(0)=0  \ \text{and }\  h''(u)> 0 \ \ \text{for} \  u>0.
\end{array}
 \]
Similarly, concave function $h$ belongs to $\mathbb{Y}\subset C^2(\mathbb{R})$ if we have
\[
\begin{array}{ccc}
\text{either} & h(0)>0  \  \text{and }\  h''(u)\le 0 \ \ \text{for} \  u>0, &   \\
 \text{or} & \ h(0)=0  \ \text{and }\  h''(u)< 0 \ \ \text{for} \  u>0.
\end{array}
 \]
For the following single equation with the Robin boundary condition,
\begin{eqnarray*}
\label{?}
\hbox{$(N_{f})$}\hskip 50pt \left\{ \begin{array}{lcl}
\hfill -\Delta u =f(u)     &\text{in}  \ \Omega, &  \\
\hfill au+b\partial_\nu u=0    &\text{on} \  \partial \Omega, &
\end{array}\right.
    \end{eqnarray*}
Brown and Shivaji in \cite{bs} proved the following result.
\begin{thm}
Suppose $f\in\mathbb{X}$ and $f'(u)>0 $ for $u>0$, then all $C^2(\Omega)$ positive solutions of $ (N_{f})$ are unstable.
\end{thm}
In the following theorem, Tertikas \cite{te} improved the above result by removing the  monotonicity condition $f'(u)>0$, for $u>0$.
\begin{thm}
(i) \,If $f\in\mathbb{X}$, then every $C^2(\Omega)$ positive solution of  $ (N_{f})$ is unstable.\\
(ii)\, If $f\in\mathbb{Y}$, then all $C^2(\Omega)$ positive solutions  of $(N_{f})$ are stable.
\end{thm}
This theorem is sharp. Since for the Gelfand non-linearity $f(u)=e^u$ which we have $f(0)>0$, by Montenegro's paper \cite{mo} we know that minimal solutions are stable.  Now, we try to prove a counterpart of these theorems for systems.
%We should mention that, recently, Montenegro \cite{mo} showed minimal solution of some elliptic systems with $f(0)>0, g(0)>0$ and $ f'',g''>0$, like Gelfand system $(N_{1,\lambda,\gamma})$,  are stable.

\begin{thm}
(i) Let $f,g\in\mathbb{X}$, $\theta\in C^1(\Omega)$ and $\lambda,\gamma\in \mathbb{R^+}$. If $(u,v) $  is a $C^2(\Omega)$ positive solution of $(N_{f,g})$, then $(u,v)$ is unstable.
\\
(ii) Let $f,g\in \mathbb{Y}$, $\theta\in C^1(\Omega)$ and $\lambda,\gamma\in\mathbb{R^+}$. Then, all $C^2(\Omega)$ positive solutions of $(N_{\lambda,\gamma})$ are stable.
\end{thm}

%\begin{lemma}
%(i) For $f\in\mathbb{X}$ we have $$f'(u)u-f(u)>0\ \text{for}\ \ u>0$$
%(ii) For $f\in\mathbb{Y}$ we have $$f'(u)u-f(u)<0\ \text{for}\ \ u>0$$
%\end{lemma}
%Proof: Just take a look at the derivative of $F(u)=f'(u)u-f(u)$. We see $F'(u)=f''(u)u$. So, for (i) we have $F(0)>0$ and $ F'(u)\ge 0,u>0 $ or  $F(0)=0$ and $ F'(u)>0,u>0 $ which says $F(u)>0,u>0$.  (ii) is similar.

\noindent{\it Proof.}\, Let $(u,v)$ be a positive semi-stable solution of $(N_{f,g})$, so there exists positive eigenfunction pair $(\phi,\psi)$ corresponding to the first eigenvalue $\eta$  such that
\begin{eqnarray*}
\label{?}
\hbox{$(S_{f,g})$}\hskip 50pt \left\{ \begin{array}{lcl}
\hfill -div(e^\theta \nabla \phi)&=&\lambda e^\theta f'(v)\psi+\eta\phi  \ \ \ \text{in} \  \Omega,   \\
\hfill -div(e^\theta \nabla \psi)&=&\gamma e^\theta g'(u)\phi+\eta\psi  \ \ \ \text{in} \  \Omega,
\end{array}\right.
  \end{eqnarray*}
with Robin  boundary conditions:
  \begin{eqnarray*}
 \left\{ \begin{array}{lcl}
\hfill a\phi+b\partial_\nu \phi&=&0  \ \ \ \text{on} \ \partial\Omega, \\
\hfill  a\psi +b\partial_\nu \psi &=&0 \ \ \ \text{on} \ \partial\Omega. \\
\end{array}\right.
  \end{eqnarray*}
Multiply the first equation of $(S_{f,g})$ by $v$ and the second equation of $(N_{f,g})$ by $\phi$ and integrate to get:
\begin{eqnarray*}
\lambda\int_{\Omega} e^\theta f'(v)v\psi+\int_{\Omega}\eta v\phi& =&-\int_{\Omega}div(e^\theta \nabla\phi) v\\
&=&-\int_{\Omega}div(e^\theta \nabla v) \phi+\int_{\partial\Omega}e^\theta\left( {\phi\partial_\nu v-v\partial_\nu \phi}\right)\\
&=&\gamma\int_{\Omega}{e^\theta g(u)\phi}+\int_{\partial\Omega} e^\theta \left({\phi\partial_\nu v-v\partial_\nu\phi}\right)\\
&=& \gamma\int_{\Omega}{ e^\theta g(u)\phi},
\end{eqnarray*}
the last equality holds, since $v$ and $\phi$ have the same boundary condition. Similarly,
\begin{eqnarray*}
\gamma \int_{\Omega} e^\theta g'(u)u\phi+\int_{\Omega}\eta u\psi& =&-\int_{\Omega}div(e^\theta \nabla\psi) u\\
&=&-\int_{\Omega}div(e^\theta \nabla u) \psi+\int_{\partial\Omega}e^\theta\left( {\psi\partial_\nu u-u\partial_\nu \psi}\right)\\
&=&\lambda\int_{\Omega}{e^\theta f(v)\psi}+\int_{\partial\Omega} e^\theta \left({\psi\partial_\nu u-u\partial_\nu\psi}\right)\\
&=& \lambda\int_{\Omega}{ e^\theta f(v)\psi}.
\end{eqnarray*}
For (i),  since $f$ is in $\mathbb{X}$, we have $f'(v)v>f(v) $. Therefore,
\begin{eqnarray*}
 \gamma\int_{\Omega}{  e^\theta g(u)\phi}&\ge&  \lambda \int_{\Omega}e^\theta f(v)\psi+\int_{\Omega}\eta v\phi\\
 &=&\gamma \int_{\Omega} e^\theta g'(u)u\phi+\eta\int_{\Omega} (u\psi+ v\phi).
\end{eqnarray*}
Suppose $(u,v)$ is a semi-stable solution, so $\eta\ge0$ and we have
\begin{equation} \label{1}
 \gamma\int_{\Omega}{  e^\theta   (g(u)-g'(u)u)\phi}\ge 0,
 \end{equation}
 and this is a contradiction, since  $g\in\mathbb{X}$, i.e., $g(u)-g'(u)u<0$.

For (ii),  since  $f$ and $g$ are in $\mathbb{Y}$, we have $f(v)>f'(v)v $ and $g(u)>g'(u)u$. Therefore,
\begin{eqnarray*}
 \gamma\int_{\Omega}{ e^\theta g(u)\phi}&\le&  \lambda \int_{\Omega} e^\theta f(v)\psi+\int_{\Omega}\eta v\phi\\
 &=&\gamma \int_{\Omega} e^\theta g'(u)u\phi+\eta\int_{\Omega} (u\psi+ v\phi),
\end{eqnarray*}
and,
\begin{eqnarray*}
 \eta\int_{\Omega} (u\psi+ v\phi)\ge \gamma\int_{\Omega}{e^\theta (g(u)-g'(u)u)\phi}>0,
\end{eqnarray*}
we conclude all eigenvalues are positive.

\hfill $\Box$

\noindent{\bf Remarks}
\begin{enumerate}

\item There is no condition on $f'$ and $g'$.

\item The same result holds for $f(u,v,x)$ and $g(u,v,x)$.

\item This theorem looks sharp. Since for the Gelfand non-linearity $f=g=\exp(.)$ which we have $f(0)>0$ and $g(0)>0$, by Montenegro's paper \cite{mo} we know that minimal solutions are stable. Moreover, he actually proved some properties of  the curve $\Gamma=\partial\Lambda$, where
 $$\Lambda :=\{(\lambda, \gamma ) \in\mathbb{R^{+}} : \text{there exists a smooth solution $(u, v)$ of}  \ (N_{f,g}) \}. $$
 We can easily get an upper bound for the curve  $ \Gamma  \subseteq \{   (\lambda,\gamma) \in\omega ; \ \lambda\gamma\le\frac{\lambda_1^2}{f'(0)g'(0)}\}$. Since $f'',g''>0$, we have $f(v)\ge f'(0)v$ and $g(u)\ge g'(0)u$. Multiply both equations  with the first eigenfunction of the operator $-\Delta$, i.e. $\phi_1$, and do integration by parts to get
$$ \lambda_1\int_{\Omega}{ \phi_1 u}\ge \lambda f'(0)\int_{\Omega}{  \phi_1  v}  \ \ \text{and} \ \   \lambda_1\int_{\Omega}{ \phi_1 v}\ge \gamma g'(0)\int_{\Omega}{  \phi_1 u}.$$

\end{enumerate}

\vspace{2mm}For the remainder of the paper, we will focus on Pohozaev type theorems for the following weighted elliptic system, in the absence of stability,
\begin{eqnarray*}
\label{?}
\hbox{$(N_{f,g})$}\hskip 50pt \left\{ \begin{array}{lcl}
\hfill -\Delta u&=& |x|^\alpha \ f(v) \ \  \text{in} \  \Omega,   \\
\hfill -\Delta v&=&  |x|^\beta  \ g(u)   \ \ \text{in} \  \Omega,
\end{array}\right.
  \end{eqnarray*}
where  $f,g\in C^1(\mathbb{R})$ and $\alpha,\beta\ge 0$, with Dirichlet boundary conditions $u=v=0   \ \text{on} \  \partial\Omega$.
\begin{lemma} \label{poho}
All non-negative $C^2(\Omega)$ solutions of  $(N_{f,g})$ satisfy
\begin{eqnarray*}
&\int_{\Omega}& (  (N+\alpha) F(v) -(N-2)\tau \ v f(v)  ) |x|^\alpha+ \\
& \int_{\Omega}& (   (N+\beta) G(u) -(N-2) (1-\tau) \  u g(u) ) |x|^\beta =\\
& \int_{\partial\Omega}& (x\cdot\nu) \partial_\nu u\partial_\nu v dS,
\end{eqnarray*}
for any $\tau$ in $[0,1]$. The functions $F$ and $G$ stand for anti-derivative functions of $f$ and $g$, respectively, i.e., $F(t):=\int_{0}^{t} f(s) ds$.
\end{lemma}
\noindent{\it Proof.} \, First, we show the following identity holds for any $C^2(\Omega)$ functions $u $ and $v$ which are zero on the boundary of $\Omega$.
\begin{eqnarray}\label{pohoiden}
\int_{\Omega} (x\cdot\nabla v) \Delta u+ (x\cdot\nabla u) \Delta v -(N-2) \nabla u\cdot\nabla v =  \int_{\partial\Omega} (x\cdot\nu) \partial_\nu u\partial_\nu v dS.
\end{eqnarray}
By a straightforward calculation, one can observe the following identities:
\begin{eqnarray*}
div (  (x\cdot\nabla  v) \nabla u &+& (x\cdot\nabla u)\nabla v ) \\
&=& (x\cdot\nabla v) \Delta u+ (x\cdot\nabla u) \Delta v +x\cdot \nabla (\nabla u\cdot\nabla v) +2 \nabla u\cdot\nabla v,
\\
div ( x(\nabla u\cdot\nabla v)) &=& N \nabla u\cdot \nabla v +x\cdot \nabla (\nabla u\cdot\nabla v).
\end{eqnarray*}
Subtract these equalities to get
\begin{eqnarray*}
div (  (x\cdot\nabla  v) \nabla u &+& (x\cdot\nabla u)\nabla v - x(\nabla u\cdot\nabla v)) \\
& =& (x\cdot\nabla v) \Delta u+ (x\cdot\nabla u) \Delta v -(N-2) \nabla u\cdot\nabla v.
\end{eqnarray*}
Integrating over $\Omega$ and using Green's formula, we get
\begin{eqnarray*}
&\int_{\Omega}&  \left(   (x\cdot\nabla  v)  \Delta u  +  (x\cdot\nabla  u) \Delta v -(N-2) \nabla u\cdot\nabla v \right) dx  =\\ &\int_{\partial\Omega} &\left(   (x\cdot\nabla  v) \nabla u + (x\cdot\nabla u)\nabla v - x(\nabla u\cdot\nabla v)  \right) \cdot \nu dS.
\end{eqnarray*}
Since $u=v=0$ on $\partial\Omega$, we have $\nabla u=\partial_\nu u \ \nu$ and $\nabla v=\partial_\nu u\ \nu$.  This proves (\ref{pohoiden}).

Now, assume $(u,v)$ is a solution of $(N_{f,g})$. Then,
\begin{eqnarray*}
 (x\cdot\nabla  v)  \Delta u + (x\cdot\nabla  u) \Delta v  &=& - (x\cdot\nabla  v) |x|^\alpha f(v)- (x\cdot\nabla  u) |x|^\beta g(u)\\
 &=&- x\cdot \nabla F(v) |x|^\alpha - x\cdot \nabla G(u) |x|^\beta\\
 &=& - div \left( x |x|^\alpha F(v)+x |x|^\beta G(u)     \right) \\
 &+& (N+\alpha) |x|^\alpha F(v) + (N+\beta) |x|^\beta G(u).
\end{eqnarray*}
Therefore, we have
\begin{equation} \label{pohoiden2} \begin{array}{ll}
\int_{\Omega}\left( (x\cdot\nabla  v)  \Delta u + (x\cdot\nabla  u) \Delta v \right)dx = \int_{\Omega} (N+\alpha) |x|^\alpha F(v)+(N+\beta) |x|^\beta G(u).
\end{array}
\end{equation}

On the other hand, by multiplying equations of $(N_{f,g})$ with $u$ and $v$, we get
\begin{eqnarray*}
\int_{\Omega} \nabla u \cdot\nabla v &=&- \int _{\Omega} v \Delta u dx= \int_{\Omega } |x|^\alpha v f(v) dx,\\
\int_{\Omega} \nabla u \cdot\nabla v &=&- \int _{\Omega} u \Delta v dx= \int_{\Omega } |x|^\beta u g(u) dx.
\end{eqnarray*}
Therefore, for any $\tau\in [0,1]$ we have
\begin{eqnarray}   \label{pohoiden3}
\int_{\Omega} \nabla u \cdot\nabla v = \int_{\Omega} \left (  \tau  |x|^\alpha v f(v) +(1-\tau)  |x|^\beta u g(u) \right ) dx.
\end{eqnarray}
Substitute (\ref{pohoiden2}) and (\ref{pohoiden3}) into (\ref{pohoiden}) to get the desired result.

$\hfill \blacksquare$ %\hfill $\Box$

%We call the function $h$ to be $\lambda$-superlinear if there exists a constant  $\lambda$ such that the %following condition satisfies:
%\begin{eqnarray*}
   % H(t):=\int_{0}^{t} h(s) ds \ge\lambda t \ h(t)\  \ \text{for all} \ t> 0, \ \text{and} \ \ h(0)=0.
%\end{eqnarray*}

\begin{thm} \label{pohozaev}
Let  $N\ge3$ and $\Omega\subset\mathbb{R}^N$ be a star-shaped and bounded domain. Then, there is no $C^2(\Omega)$ positive solution for $(N_{f,g})$, provided $f$ and $g$ satisfy the following conditions for all $t>0$:
\begin{eqnarray*}
  (N+\alpha) F(t) &\le& a_1 t f(t), \\
   (N+\beta) G(t) &\le& a_2 t g(t),
\end{eqnarray*}
where $a_1,a_2\in \mathbb{R}^+$ and $a_1+a_2=N-2$.

\end{thm}

\noindent{\it Proof.}\,  Since $\Omega$ is a star-shaped domain, we have $\int_{\partial\Omega} (x\cdot\nu) \partial_\nu u\partial_\nu v dS> 0$. Therefore, Lemma \ref{poho} leads us to a contradiction if we take
 $  (N+\alpha) F(v) -(N-2)\tau \ v f(v)\le 0 $  and  $(N+\beta) G(u) -(N-2) (1-\tau) \  u g(u)\le 0$.

$\hfill \blacksquare$ %\hfill $\Box$

Therefore, for the weighted Lane-Emden system we have the following result.
\begin{cor} Let $f(v)=v^p$ and $g(u)=u^q$ for $p,q\ge 1$ in  Theorem \ref{pohozaev}. Then, there is no positive classical solution for $(N_{f,g})$ in dimensions
$$\frac{N+\alpha}{p+1}+\frac{N+\beta}{q+1}\le N-2. $$

\end{cor}

\noindent{\bf Remark} \,  For the fourth order equation, i.e.  $\alpha=0, p=1$ and $N>4$, there is no positive solution for $q \ge \frac{N+4+2\beta}{N-4}$ and also for the weighted Lane-Emden equation, i.e., $p=q$ and $\alpha=\beta$, there is no positive  solution for $q \ge \frac{N+2+2\beta}{N-2}$. This equation is also called the H\'{e}non equation, see \cite{n}.\\

\noindent{\bf Acknowledgement}. I would like to thank Professor Nassif Ghoussoub, my supervisor, for his constant support. Also, many thanks to friend of mine Dr. Craig Cowan for his comments.


\small
\begin{thebibliography}{00}

\bibitem{bs}  K. J. Brown, R. Shivaji;
 \emph{Instability of nonnegative solutions for a class of semipositone problems}, Proc. Amer. Math. Soc. {\bf112} (1991), 121-124.

\bibitem{cw} C. Cowan;
\emph{Regularity of the extremal solutions in a Gelfand system problem}, Adv. Nonlinear Stud. {\bf 11} (2011), 695-700.

 \bibitem{cen}  C. Cowan,  P. Esposito, N. Ghoussoub;
 \emph{Regularity of extremal solutions in fourth order nonlinear eigenvalue problems on general domains}, Discrete and Continuous Dynamical Systems - Series A, {\bf 28} (2010), 1033-1050.



\bibitem{e1} P. Esposito, \emph{Linear instability of entire solutions for a class of non-autonomous elliptic equations,}  Proc. Roy. Soc. Edinburgh Sect. A,  {\bf 138}  (2008), 1005-1018.


\bibitem{eggu1}
 P. Esposito, N. Ghoussoub, Y. Guo;
\emph{Compactness along the branch of semistable and unstable solutions for an elliptic problem with a singular nonlinearity}, Comm. Pure Appl. Math. {\bf 60} (2007), 1731-1768.

\bibitem{eggu2} P. Esposito, N. Ghoussoub, Y. Guo;
 {Mathematical Analysis of Partial Differential Equations Modeling Electrostatic MEMS}, Courant Lecture Notes in Mathematics {\bf 20}, Courant Institute of Mathematical Sciences, New York; American Mathematical Society, Providence, RI, 2010, xiv+318 pp.

\bibitem{f} A. Farina;
 \emph{On the classification of solutions of the Lane-Emden equation on unbounded domains of $\Bbb R^N$}, J. Math. Pures Appl. (9) {\bf 87} (2007), 537-561.

\bibitem{gs} B. Gidas,  J. Spruck,
\emph{Global and local behavior of positive solutions of nonlinear
elliptic equations}, Commun. Pure Appl. Math. {\bf 34} (1981), 525-598.

\bibitem{m} E. Mitidieri,
\emph{Nonexistence of positive solutions of semilinear elliptic systems in $\mathbb{R}^{N}$}, Differential Integral Equations {\bf 9} (1996), 465-479.

\bibitem {mp} E. Mitidieri, S. I. Pokhozhaev;
\emph{A priori estimates and the absence of solutions of nonlinear partial
differential equations and inequalities}, Tr. Mat. Inst. Steklova, (2001) 234:1-384.


\bibitem{mo} M. Montenegro;
 \emph{Minimal solutions for a class of elliptic systems}, Bull. London Math. Soc. {\bf 37} (2005), 405-416.


\bibitem{n} W. M. Ni,
\emph{A nonlinear Dirichlet problem on the unit ball and its applications,}
Indiana Univ. Math. J. {\bf 31} (1982), 801-807.

 \bibitem{pqs} P. Pol\'acik, P. Quittner,  Ph. Souplet,
 \emph{Singularity and decay estimates in superlinear problems via Liouville-type theorems, Part I: elliptic systems}, Duke Math. J. {\bf 139} (2007), 555-579.


   \bibitem{qs} P. Quittner,  Ph.  Souplet;
 {Superlinear Parabolic Problems. Blow-Up, Global Existence and Steady
States}, Birkhauser Verlag, Basel, 2007.


  \bibitem{so} Ph. Souplet;
 \emph{The proof of the Lane-Emden conjecture in four space dimensions},  Adv. Math. {\bf 221} (2009), 1409-1427.


 \bibitem{sz} J. Serrin, H. Zou;
 \emph{Non-existence of positive solutions of Lane-Emden systems}, Differential Integral Equations {\bf 9} (1996), 635-653.

\bibitem{te} A. Tertikas;
 \emph{Stability and Instability of Positive Solutions of Semipositone Problems}, Proc. Amer. Math. Soc.  {\bf 114} (1992), 1035-1040.

\bibitem{wy} J. Wei, D. Ye;
 \emph{Liouville theorems for finite Morse index solutions of biharmonic problem}, preprint.


\end{thebibliography}
\end{document}